\RequirePackage{fix-cm}
\documentclass[smallextended]{svjour3}       % onecolumn (second format)
\smartqed  % flush right qed marks, e.g. at end of proof
\usepackage{graphicx}
\begin{document}

\title{ the Transformation of Sieve Function%\thanks{Grants or other notes
}


\author{Jinzhu Han%         \and
        %Second Author %etc.
}


\institute{J. Han \at
              Room 1611, Wang Kezhen Building, No.52, Haidian District, Beijing 100080, People's Republic of China \\
              Tel.: +86-13314863845\\
              \email{han87654321@sina.com}           %  \\
           \and
           %S. Author \at
              %second address
}

\date{Received: date / Accepted: date}

\maketitle

\begin{abstract}
 In this paper, we introduced the theory of the sieve function transformation. Using the principle of sieve function transformation, we improved sieve method, and obtained the difference range of similar sieve function values. For this, we improved the prime theorem for arithmetic series. The theory of sieve function transformation can also be used to research many problems in number theory.
\keywords{Sieve method\and Transformation of Sieve Function \and Prime Number Theorem}
\subclass{11M26, 11N35,11N36,11P32}
\end{abstract}

\section{Introduction}
\label{intro}
  The sieve method is a basic method of number theory.The sieve method is commonly used to study many number theory problems such as the prime number theorem, Goldbach's conjecture, twin prime number conjecture, etc. $^{[1-9]}$.

 In this paper, we will introduce  the transformation of sieve function to improve sieve method. Many notations will be used in this paper.
 
  \textbf{Notations}

   $B\ll A$, $B=O(A)$: there be positive constant $c$ such that $|B|\leq cA$, for $A>0$
   
   $\mu(d)$ : $ M\ddot{o}bius$ function
   
   $\chi(n)$ : Dirichlet characteristics
   
    $\nu(d)$ : the number of different prime divisors of $d$
    
    $\phi(q)$ : Euler function
    
    $\forall\lambda_{i}$ : all of  $\lambda_{i}$

\section{  The Transformation of Sieve Function }
In number theory, let there exist a integer set  $ A=\left \{ a:a\leq x, \right \} $, $ A_{d}=\left \{ a:a|d,a\in A \right \} $, then the  sieve function be general represented as
\begin{equation}
S(A;P(z),z)=\sum _{a\in A}\sum _{d|(a,P(z))}\mu (d)=\sum _{d|(P(z))}\mu (d)\left | A_{d} \right |
\label{eq}
\end{equation} 
where $ |A_{d}| $ is size of set $ A_{d} $, $\mu(d)$ is $ M\ddot{o}bius$ function, 
\begin{equation}
P(z)=\prod _{p\leq z}p
\label{eq}
\end{equation}

  However,  It is difficult to simple calculate  the value of sieve function $S(A;P(z),z)$ in the general case. In order to improve the sieve method, we will introduce transformation of sieve function.

The sieve function can be written as another form
\begin{equation}
 S(A;P(z),z)=|A|+\sum _{1<d|(P(z))}\mu (d)\left | A_{d} \right |
\label{eq}
\end{equation}
let us introduce an integer $\lambda_{i}$, that be called as transformation factor, such that 
\begin{equation}
|A=|A'|
\label{eq}
\end{equation}
and
\begin{equation}
|A_{d}|\to |A'_{d}|= |A_{d}|+\lambda_{i},1<d|q_{i}
\label{eq}
\end{equation}
then it have
\begin{equation}
\sum _{1<d|q_{i}}\mu (d)|A'_{d}|=\sum _{1<d|q_{i}}\mu (d)(|A_{d}|+\lambda_{i})=\sum _{1<d|q_{i}}\mu (d)|A_{d}|+\sum _{1<d|q_{i}}\mu (d)\lambda_{i}
\end{equation}
because
\begin{equation}
\sum _{d|q_{i}}\mu (d)\lambda_{i}=0
\end{equation}
\begin{equation}
\sum _{1<d|q_{i}}\mu (d)\lambda_{i}=-\lambda_{i}
\end{equation}
so it have
\begin{equation}
S(A';P(z),z)=S(A;P(z),z)+\sum^{n}_{i=1}\sum _{1<d|q_{i}}\mu (d)\lambda_{i}=S(A;P(z),z)-\sum^{n}_{i=1} \lambda_{i}
\label{eq}
\end{equation}
that be called as transformation of sieve function in this paper, and written as
\begin{equation}
S(A;P(z),z)\to S(A';P(z),z)
\label{eq}
\end{equation}
The sieve function transformation is essentially a numerical transformation of the sieve function, and there are some simple properties as follows:

Property 1, reversibility of sieve function transformation:

If there is a sieve function transformation,
\begin{equation}
S(A;P(z),z)\to  S(A';P(z),z)
\label{eq}
\end{equation}
then there must be an inverse transformation,
\begin{equation}
S(A';P(z),z)\to  S(A;P(z),z)
\label{eq}
\end{equation}
It can be defined as
\begin{equation}
|A'|=|A|
\label{eq}
\end{equation}
and
\begin{equation}
|A'_{d}|\to |A_{d}|= |A'_{d}|-\lambda_{i},1<d|q_{i}
\label{eq}
\end{equation}
where $q_{i}|P(z)$ ,so we have
\begin{equation}
S(A;P(z),z)=S(A';P(z),z)-\sum^{n}_{i=1}\sum _{1<d|q_{i}}\mu (d)\lambda_{i}=S(A';P(z),z)+\sum^{n}_{i=1} \lambda_{i}
\label{eq}
\end{equation}
The sieve function of mutual transformation, also known as the similarity sieve function in this article, is written as
\begin{equation}
S(A;P(z),z)\sim S(A';P(z),z)
\label{eq}
\end{equation}

Property 2, transitivity of sieve function transformation:

If there is a sieve function transformation,
\begin{equation}
S(A;P(z),z)\to  S(A';P(z),z)
\label{eq}
\end{equation}
and
\begin{equation}
S(A';P(z),z)\to  S(A'';P(z),z)
\label{eq}
\end{equation}
then there must be a sieve function transformation
\begin{equation}
S(A;P(z),z)\to  S(A'';P(z),z)
\label{eq}
\end{equation}
namely,for $S(A;P(z),z)\sim  S(A';P(z),z)$,$S(A';P(z),z)\sim  S(A'';P(z),z)$,then it have $S(A;P(z),z)\sim  S(A'';P(z),z)$.

Property 3: Identity of sieve function transformation:

The simplest sieve function transformation is as follows,
\begin{equation}
S(A;P(z),z)\to|A|
\label{eq}
\end{equation}
and
\begin{equation}
|A'|\to S(A';P(z),z)
\label{eq}
\end{equation}
so let $|A|=|A'|$,then there must be at least one way to screen the function transformation $S(A;P(z),z)\to  S(A';P(z),z)$;namely,let $|A|=|A'|$,then must have $S(A;P(z),z)\sim  S(A';P(z),z)$.

Conversely,let $S(A;P(z),z)\sim  S(A';P(z),z)$,then must have $|A|=|A'|$.

This is the identity of the sieve function transformation.

\section {Concepts related to sieve function transformation}

For the convenience of description, for $S(A;P(z),z)\sim  S(A';P(z),z)$, the transformation relationship can also be written as
\begin{equation}
|A'_{d}|-|A_{d}|=\lambda_{i},1<d|q_{i}
\label{eq}
\end{equation}
transformation factor $\lambda_{i}$ can also be written as $\lambda(q_{i})$,
\begin{equation}
|A'_{d}|-|A_{d}|=\lambda(q_{i}),1<d|q_{i}
\label{eq}
\end{equation}
The transformation value depends on the transformation factor $\lambda(q_{i})$,we also refer to the transformation factor $\lambda(q_{i})$ as the difference factor;$\nu(q_{i})$ represents the number of different prime factors contained in $q_{i}$. We also refer to $\nu (q_ {i}) $as the dimension of the differential factor $\lambda(q_{i})$.

The  sieve function can be transformed as many times, that will form a group of transformations. There are some useful concepts about sieve transformation in this paper.

 \textbf{1.} Standard transformation of sieve function: let there be transformation $S(A;P(z),z)\to  S(A';P(z),z)$
\begin{equation}
S(A';P(z),z)=S(A;P(z),z)+\sum^{n}_{i=1}\sum _{1<d|q_{i}}\mu (d)\lambda_{i}=S(A;P(z),z)-\sum^{n}_{i=1} \lambda_{i}
\label{eq}
\end{equation}
such that the signs of  transformation factors $\lambda_{i}$ be all same, namely
\begin{equation}
 \forall\lambda_{i}\leq 0
\label{eq}
\end{equation}
or 
\begin{equation}
\forall \lambda_{i}\leq 0
\label{eq}
\end{equation}
then which transformation be called as standard transformation of sieve function. 

\textbf{2.} Identical transformation of sieve function: let there be transformation $S(A;P(z),z)\to S(A';P(z),z)$,
\begin{equation}
S(A';P(z),z)=S(A;P(z),z)+\sum^{n}_{i=1}\sum _{1<d|q_{i}}\mu (d)\lambda_{i}=S(A;P(z),z)-\sum^{n}_{i=1} \lambda_{i}
\label{eq}
\end{equation}
such that 
\begin{equation}
\sum^{n}_{i=1} \lambda_{i}=0
\end{equation}
so it have
\begin{equation}
S(A';P(z),z)=S(A;P(z),z)
\label{eq}
\end{equation}
the value of sieve function $S(A;P(z),z)$ be constant, which be called as identical transformation of sieve function, and  be written as $S(A';P(z),z)=S(A;P(z),z)$.
 
  \textbf{3.} Special transformation of sieve function: let there be transformation $S(A;P(z),z)\to S(A';P(z),z)$, such that satisfying two conditions
\begin{equation}
S(A;P(z),z)=S(A';P(z),z)
\label{eq}
\end{equation}
\begin{equation}
\sum _{p\leq z}|A_{p}|=\sum _{p\leq z}|A'_{p}|
\label{eq}
\end{equation}
then transformation  $S(A;P(z),z)\to S(A';P(z),z)$ be called as special transformation of sieve function. There are two basic ways  of  special transformations of sieve functions as follows.

 \textbf{The first way:} The dimensions of the transformation factor and the difference factor are the same.
 
 For example:let us  make transformation $S(A;P(z),z)\to S(A';P(z),z)$, such that 
\begin{equation}
|A_{d}|\to |A'_{d}|= |A_{d}|+\lambda,1<d|q_{1} 
\end{equation}
\begin{equation}
|A_{d}|\to |A'_{d}|= |A_{d}|-\lambda,1<d|q_{2} 
\end{equation}
and $\nu(q_{1})=\nu(q_{2})$ ,then we have
\begin{equation}
S(A;P(z),z)=S(A';P(z),z)
\label{eq}
\end{equation}
\begin{equation}
\sum _{p\leq z}|A'_{p}|=\sum _{p\leq z}|A_{p}|+\lambda(\nu(q_{1})-\nu(q_{2}))
\label{eq}
\end{equation}
so,the sieve function transformation $S(A';P(z),z)$ is a special identity sieve function transformation. In this type of special identity transformation,the dimension of the transformation factor and the difference factor are the same $\nu(q_{1}=\nu(q_{2}))$.

\textbf{The second way:}  the dimensions of transformation factor and difference factor are not the same.

For example, if we take A transformation $S(A; P(z),z)\to S(A'; P(z),z)$, so
\begin{equation}
|A_{d}|\to |A'_{d}|= |A_{d}|-\lambda,1<d|q_{1}
\end{equation}
\begin{equation}
|A_{d}|\to |A'_{d}|= |A_{d}|+\lambda,1<d|q_{2}
\end{equation}
\begin{equation}
|A_{d}|\to |A'_{d}|= |A_{d}|+\lambda,1<d|q_{3}
\end{equation}
\begin{equation}
|A_{d}|\to |A'_{d}|= |A_{d}|-\lambda,1<d|q_{4}
\end{equation}
Where $\nu (q_ {1}) = \nu (q_ {2}) + 1 $, $\nu (q_ {3}) = \nu (q_ {4}) + 1 $, satisfy the following equation
\begin{equation}
\lambda(\nu(q_{1})-\nu(q_{2}))-\lambda(\nu(q_{3})-\nu(q_{4}))=0
\label{eq}
\end{equation}
Then we have
\begin{equation}
S(A; P(z),z)=S(A'; P(z),z)
\label{eq}
\end{equation}
\begin{equation}
\sum _{p\leq z}|A'_{p}|=\sum _{p\leq z}|A_{p}|-\lambda(\nu(q_{1})-\nu(q_{2}))+\lambda(\nu(q_{3})-\nu(q_{4}))=\sum _{p\leq z}|A_{p}|
\label{eq}
\end{equation}
Thus, the sieve function transformation $S(A; P(z),z)=S(A'; P(z),z)$is also a special identity sieve function transformation. In this type of special identity transformation, the dimensions of the transformation factor and the difference factor are not the same.

The special special identity transformation of the sieve function is a very important class of transformations, using the first type of special sieve function transformation, we can merge the difference factors that have the same dimension horizontally. For example, if $S(A; P(z),z)\to S(B; P(z),z)$has difference factors of the same dimension
\begin{equation}
|A_{d}|- |B_{d}|= \lambda_{1},1<d|q_{1}
\end{equation}
\begin{equation}
|A_{d}|- |B_{d}|= \lambda_{2},1<d|q_{2}
\end{equation}
where $\nu(q_{1})=\nu(q_{2})$,we can choose $S(A; P(z),z)\to S(A'; P(z),z)$,such that
\begin{equation}
|A_{d}|\to |A'_{d}|= |A_{d}|-\lambda_{1},1<d|q_{1}
\end{equation}
\begin{equation}
|A_{d}|\to |A'_{d}|= |A_{d}|+\lambda_{1},1<d|q_{2}
\end{equation}
Then we have
\begin{equation}
S(A; P(z),z)=S(A'; P(z),z)
\label{eq}
\end{equation}
\begin{equation}
\sum _{p\leq z}|A'_{p}|=\sum _{p\leq z}|A_{p}|-\lambda_{1}(\nu(q_{1})-\nu(q_{2}))=\sum _{p\leq z}|A_{p}|
\label{eq}
\end{equation}
Thus, the sieve function transformation $S(A; P(z),z)=S(A'; P(z),z)$ is a special identity sieve function transformation, and the sieve function transformation $S(A'; P(z),z)\to S(B;P(z),z)$ becomes
\begin{equation}
|A'_{d}|- |B_{d}|= 0,1<d|q_{1}
\end{equation}
\begin{equation}
|A'_{d}|- |B_{d}|= \lambda_{2}+\lambda_{1},1<d|q_{2}
\end{equation}
In this way, the difference factors of the same dimension are merged into one of them, which can be called horizontal merging of the difference factors.

Some difference factors with the same dimension may be positive or negative. Through the horizontal merger  of difference factors, the difference factors with different symbols can be offset each other, so that all difference factors with the same dimension are greater than or equal to zero, or all are less than or equal to zero. Through the horizontal merging of the difference factors, the difference factors with the same dimension can also be averaged, so that the difference factors with the same dimension are roughly the same.

Using the second type of special sieve function transformation, we can vertical merge differential factors with different dimensions. For example:

If $S(A; P(z),z)\to S(B; P(z),z)$has many different factors
\begin{equation}
|A_{d}|- |B_{d}|= \lambda_{1},1<d|q_{1}
\end{equation}
\begin{equation}
|A_{d}|- |B_{d}|= \lambda_{2},1<d|q_{2}
\end{equation}
\begin{equation}
|A_{d}|- |B_{d}|= \lambda_{3},1<d|q_{3}
\end{equation}
\begin{equation}
|A_{d}|- |B_{d}|= \lambda_{4},1<d|q_{4}
\end{equation}
Where $\nu (q_ {1}) = 5 $, $\nu (q_ {2}) = 4 $, $\nu (q_ {3}) = 2 $, $\nu (q_ {4}) = 1 $, we can select screen function transform $S (A; P(z),z)\to S(A'; P(z),z)$, so
\begin{equation}
|A_{d}|\to |A'_{d}|= |A_{d}|-\lambda_{1},1<d|q_{1}
\end{equation}
\begin{equation}
|A_{d}|\to |A'_{d}|= |A_{d}|+\lambda_{1},1<d|q_{2}
\end{equation}
\begin{equation}
|A_{d}|\to |A'_{d}|= |A_{d}|+\lambda_{1},1<d|q_{3}
\end{equation}
\begin{equation}
|A_{d}|\to |A'_{d}|= |A_{d}|-\lambda_{1},1<d|q_{4}
\end{equation}
Then we have
\begin{equation}
S(A; P(z),z)=S(A'; P(z),z)
\label{eq}
\end{equation}
\begin{equation}
\sum _{p\leq z}|A'_{p}|=\sum _{p\leq z}|A_{p}|-\lambda_{1}(\nu(q_{1})-\nu(q_{2}))+\lambda_{1}(\nu(q_{3})-\nu(q_{4}))=\sum _{p\leq z}|A_{p}|
\label{eq}
\end{equation}
Thus, the sieve function transformation $S(A; P(z),z)=S(A'; P(z),z)$is special identity sieve function transformation, and the difference factor of sieve function transformation $S(A'; P(z),z)\to S(B;  P(z),z)$ becomes
\begin{equation}
|A'_{d}|- |B_{d}|= 0,1<d|q_{1}
\end{equation}
\begin{equation}
|A'_{d}|- |B_{d}|= \lambda_{2}+\lambda_{1},1<d|q_{2}
\end{equation}
\begin{equation}
|A'_{d}|- |B_{d}|= \lambda_{3}+\lambda_{1},1<d|q_{3}
\end{equation}
\begin{equation}
|A'_{d}|- |B_{d}|= \lambda_{3}-\lambda_{1},1<d|q_{4}
\end{equation}
The difference factor with a larger dimension is zero, that is,$ |A'_{d}|-|B_{d}|=0,1<d|q_{1} $ . In this way, the difference factor with a larger dimension becomes zero and is merged into the difference factor with a smaller dimension, which is called vertical merger of difference factors. Using the second type of special identity sieve function transformation, the difference factors can be vertically merged many times, so that the difference factors with larger dimensions become zero, and finally only the difference factors with dimensions less than or equal to 2 can be retained, that is, only the difference factors $\lambda(q_{i})$ with $\nu(q_{i})\leq 2$ can be retained.

\section{Theorem relating to sieve function transformation}

\begin{theorem}
Let $S(A; P(z),z)\to S(B; P(z),z)$is a standard sieve function transformation, then there is
\begin{equation}
S(A; P(z),z)-S(B; P(z),z)\ll\left|\sum _{p\leq z}|A_{p}|-\sum _{p\leq z}|B_{p}|\right|
\label{eq}
\end{equation}
\end{theorem}

\begin{proof}
Let $S(A;P(z),z)\to S(B;P(z),z)$ be a standard sieve function transformation,
\begin{equation}
S(B;P(z),z)=S(A;P(z),z)+\sum^{n}_{i=1}\sum _{1<d|q_{i}}\mu (d)\lambda_{i}=S(A;P(z),z)-\sum^{n}_{i=1} \lambda_{i}
\label{eq}
\end{equation}
and all of the transformation factors, that is, the difference factor $\lambda_{i}$, have the same sign, thus
\begin{equation}
|S(A; P(z),z)-S(B; P(z),z)|=\sum^{n}_{i=1} |\lambda_{i}|
\label{eq}
\end{equation}
obviously, for $\forall\lambda_{i}\leq 0$, or $\forall\lambda_{i}\geq 0$, then it have
\begin{equation}
\sum^{n}_{i=1} |\lambda_{i}|\leq\left|\sum _{p\leq z}|A_{p}|-\sum _{p\leq z}|B_{p}|\right|
\label{eq}
\end{equation}
so for $\forall\lambda_{i}\leq 0$, or $\forall\lambda_{i}\geq 0$, then it have
\begin{equation}
S(A; P(z),z)-S(B; P(z),z)\ll\left|\sum _{p\leq z}|A_{p}|-\sum _{p\leq z}|B_{p}|\right|
\label{eq}
\end{equation}
Theorem 1 is proved.
\end{proof}

\begin{theorem}
Let $S(A; P(z),z)\sim S(B; P(z),z)$, namely, there be a sieve function transformation $S(A; P(z),z)\to S(B; P(z),z)$, then it have
\begin{equation}
S(A; P(z),z)-S(B; P(z),z)\ll \sum _{p\leq z}\left||A_{p}|-|B_{p}|\right|+R_{0}
\label{eq}
\end{equation}
where $R_{0}$ be special error term,
\begin{equation}
R_{0}\ll \sum _{p\leq z}\log |A|\ll \sum _{p\leq z}\log^{2} |A|
\label{eq}
\end{equation}

\end{theorem}

\begin{proof}
let $S(A;P(z),z)\sim S(B;P(z),z)$,there be standard sieve function transformation $S(A;P(z),z)\to S(B;P(z),z)$,then according to Theorem 1,we have 
\begin{equation}
S(A;P(z),z)-S(B;P(z),z)\ll\left|\sum _{p\leq z}|A_{p}|-\sum _{p\leq z}|B_{p}|\right|
\label{eq}
\end{equation}
since
\begin{equation}
 \left|\sum _{p\leq z}|A_{p}|-\sum _{p\leq z}|B_{p}|\right|\leq\sum _{p\leq z}\left||A_{p}|-|B_{p}|\right|
\label{eq}
\end{equation}
thus
\begin{equation}
S(A;P(z),z)-S(B;P(z),z)\ll \sum _{p\leq z}\left||A_{p}|-|B_{p}|\right|
\label{eq}
\end{equation}

Let $S(A;P(z),z)\to S(B;P(z),z)$ be not standard sieve function transformation,the difference factor can be positive or negative, so we can transform through a series of sieve functions$S(A;P(z),z)\to S(A'';P(z),z)$,such that $S(A'';P(z),z)\to S(B;P(z),z)$ become standard sieve function transformation,then according to Theorem 1, we can prove theorem 2.

There are many ways to achieve this goal. Here we choose the two types of special sieve function transformations described above.

We can combine different difference factors in multiple steps.

The first step is to vertical merge the difference factors.

For example, if $S(A; P(z),z)\to S(B; P(z),z)$is not a standard sieve function transformation, and there are many difference factors
\begin{equation}
|A_{d}|- |B_{d}|= \lambda_{1},1<d|q_{1}
\end{equation}
\begin{equation}
|A_{d}|- |B_{d}|= \lambda_{2},1<d|q_{2}
\end{equation}
\begin{equation}
\cdots\cdots
\end{equation}
\begin{equation}
|A_{d}|- |B_{d}|= \lambda_{3},1<d|q_{3}
\end{equation}
\begin{equation}
|A_{d}|- |B_{d}|= \lambda_{4},1<d|q_{4}
\end{equation}
Where $\nu(q_{1})$is the maximum value,$\nu(q_{2})=\nu(q_{1})-1$,$\nu(q_{3})=2$,$\nu(q_{4})=1$.We can choose $S(A;P(z),z)\to S(A';P(z),z)$,such that
\begin{equation}
|A_{d}|\to |A'_{d}|= |A_{d}|-\lambda_{1},1<d|q_{1}
\end{equation}
\begin{equation}
|A_{d}|\to |A'_{d}|= |A_{d}|+\lambda_{1},1<d|q_{2}
\end{equation}
\begin{equation}
|A_{d}|\to |A'_{d}|= |A_{d}|+\lambda_{1},1<d|q_{3}
\end{equation}
\begin{equation}
|A_{d}|\to |A'_{d}|= |A_{d}|-\lambda_{1},1<d|q_{4}
\end{equation}
for this, the sieve function transformation $S(A;P(z),z)=S(A';P(z),z)$ be a special identity sieve function transformation,and the difference factor of  sieve function transformation $S(A';P(z),z)\to S(B;P(z),z)$ become
\begin{equation}
|A'_{d}|- |B_{d}|= 0,1<d|q_{1}
\end{equation}
\begin{equation}
|A'_{d}|- |B_{d}|= \lambda_{2}+\lambda_{1},1<d|q_{2}
\end{equation}
\begin{equation}
\cdots\cdots
\end{equation}
\begin{equation}
|A'_{d}|- |B_{d}|= \lambda_{3}+\lambda_{1},1<d|q_{3}
\end{equation}
\begin{equation}
|A'_{d}|- |B_{d}|= \lambda_{3}-\lambda_{1},1<d|q_{4}
\end{equation}
the difference factor with larger dimensions be equal to zero,namely,$|A'_{d}|- |B_{d}|= 0,1<d|q_{1}$,that causes the difference factors with larger dimensions to become zero and merge into the difference factors with smaller dimensions.

In this way, we can use the second type of special sieve function transform, multiple vertical merge difference factors, to change the difference factor of the largest difference factor to zero, and merge to the smaller difference factors of the dimension, so many vertical merge difference factors,in the end, the difference factor of the sieve function transformation $S(A';P(z),z)\to S(B;P(z),z)$ is divided into the lower difference factor term, and finally there is only  the two types of differential factors $\lambda(p_{i}p_{j})$, and $\lambda(p_{i})$, namely
\begin{equation}
\lambda(q_{i}))=0,\nu(q_{i}))> 2
\label{eq}
\end{equation}

In the second step, the horizontal merge difference factor, we can use the first type of special  transformation$S(A';P(z),z)= S(A'';P(z),z)$,horizontal merge difference factor$\lambda(p_{i}p_{j})$,and $\lambda(p_{i})$,and the number of these difference factors is roughly equalized,for this such that the difference factor of sieve function transformation$S(A'';P(z),z)\to S(B;P(z),z)$,$\lambda(p_{i}p_{j})$ become
\begin{equation}
\forall \lambda(p_{i}p_{j})\leq 0,
\label{eq}
\end{equation}
or
\begin{equation}
\forall \lambda(p_{i}p_{j})\geq 0
\label{eq}
\end{equation}
In general cases,when $\forall\lambda(p_{i}p_{j})\leq 0$,
\begin{equation}
\forall \lambda(p_{i})\leq 0,
\label{eq}
\end{equation}
when $\forall\lambda(p_{i}p_{j})\geq 0$,
\begin{equation}
\forall \lambda(p_{i})\geq 0,
\label{eq}
\end{equation}
for this, $S(A'';P(z),z)\to S(B;P(z),z)$ be standard sieve function transformation,according to Theorem 1, we have
\begin{equation}
S(A'';P(z),z-S(B;P(z),z))\ll \left|\sum _{p\leq z}|A''_{p}|-\sum _{p\leq z}|B_{p}|\right|
\label{eq}
\end{equation}
However,the set $A$ may have the multiplication number of primes $a=p^{\alpha}$,$\alpha=2,3,4,\cdots$, sieve function transformation $S(A'';P(z),z)\to S(B;P(z),z)$ may have an exception,namely,when $\lambda(p_{i}p_{j})\leq0$ ,may have an exception$\lambda(p_{i})>0$ .We can calculate the number of multiplication number of primes in natural number set $A=\{ a:a\leq x \}$,
\begin{equation}
 \sum _{^{a=p^{\alpha}}_{p\leq z}}1\leq \sum _{p\leq z}\log |A|
\label{eq}
\end{equation}
If there is multiplication number of primes that may lead to the exception $a=p^{\alpha}$,then this special error term $R_{0}$ should be less than or equal to $\sum _{p\leq z}\log |A|$,
\begin{equation}
R_{0}\ll \sum _{p\leq z}\log |A|\ll \sum _{p\leq z}\log^{2} |A|
\label{eq}
\end{equation}
so we have
\begin{equation}
S(A'';P(z),z-S(B;P(z),z))\ll \left|\sum _{p\leq z}|A''_{p}|-\sum _{p\leq z}|B_{p}|\right|+R_{0}
\label{eq}
\end{equation}
because $S(A;P(z),z)= S(A';P(z),z)$,$S(A';P(z),z)= S(A'';P(z),z)$ be all special identical sieve function transformations.
\begin{equation}
 \sum _{p\leq z}|A_{p}|=\sum _{p\leq z}|A'_{p}|=\sum _{p\leq z}|A''_{p}|
\label{eq}
\end{equation}
so
\begin{equation}
S(A;P(z),z)-S(B;P(z),z)\ll\left|\sum _{p\leq z}|A_{p}|-\sum _{p\leq z}|B_{p}|\right|+R_{0}
\label{eq}
\end{equation}
since
\begin{equation}
 \left|\sum _{p\leq z}|A_{p}|-\sum _{p\leq z}|B_{p}|\right|\leq\sum _{p\leq z}\left||A_{p}|-|B_{p}|\right|
\label{eq}
\end{equation}
thus
\begin{equation}
S(A;P(z),z)-S(B;P(z),z)\ll\sum _{p\leq z}\left||A_{p}|-|B_{p}|\right|+R_{0}
\label{eq}
\end{equation}
Theorem 2 be proved.
\end{proof}

\section{Improvement of the prime number theorem in the arithmetic sequence}

In the number theory,the function $\pi(x;q,l)$ be defined as follows:
\begin{equation}
\pi(x;q,l)=\sum_{^{p\equiv l mod(q)}_{p\leq x}}1
\end{equation}
where $(q,l)=1$, $l<q$ ,that represents the number of primes in the arithmetic sequence.Because $\pi(x;2,1)=\pi(x)$,
\begin{equation}
\pi(x)=\sum_{p\leq x}1
\end{equation}
This paper mainly studies the situation of $q>2$. The mathematicians have proved the following theorem $^{[8]}$.

\begin{theorem}
Let $(q,l)=1$, $l<q$ ,$q\leq \log^{2}x$ ,then there be a calculable constant $c$, such that 
\begin{equation}
\pi(x;q,l)=\frac{\pi(x)}{\phi(q)}+O\left(xe^{c\sqrt{\log x}}\right)
\end{equation}
where $\phi(q)$ be Euler function,
\begin{equation}
\phi(q)=q\prod_{p|q}\left(1-\frac{1}{p}\right)
\label{eq}
\end{equation}
\end{theorem}

In this article, we apply the transformation theory of the sieve function, and improve the prime number theorem in the arithmetic sequence, and prove the following theorem. 

\begin{theorem}
Let $(q,l)=1$, $l<q$ , then it have
\begin{equation}
\pi(x;q,l)=\frac{\pi(x)}{\phi(q)}+O\left(\sqrt{x}\right)
\label{eq}
\end{equation}
\end{theorem}
\begin{proof}

As we all know,in natural numbers, the primes are roughly uniform distributed among the irreducible residue classes of $ mod(q) $.The number of the irreducible residue classes of $ mod(q) $ be equal to Euler function $\phi(q)$,for this,let $(q,l_{i})=1$, $l_{i}<q$ ,then it have
\begin{equation}
\pi(x;q,l_{i})=\frac{\pi(x)}{\phi(q)}+R_{i}
\label{eq}
\end{equation}
where $R_{i}$ be error term.Let the primes are average distributed among the irreducible residue classes of $ mod(q) $,then
\begin{equation}
\pi(x;q,l_{i})=\frac{\pi(x)}{\phi(q)}
\label{eq}
\end{equation}
namely,the mean value of $\pi(x;q,l_{i})$ be equal to $\frac{\pi(x)}{\phi(q)}$.

The function $\pi(x;q,l_{i})$ can also be expressed as a sieve function.

For example, let us  set up 

$ A^{i}=\left \{a:n=kq+l_{i}, n\leq x\right \} $, $ A^{i}_{d}=\left \{a:a\in A^{i}, a\equiv0 mod (d)\right \} $, $(q,l_{i})=1$,$l_{i}<q$, $z=\sqrt{x}$, then we have
\begin{equation}
S_{0}(A^{i};P(z),\sqrt{x})=| A^{i}|+\sum _{1<d|(P(z))}\mu (d) | A^{i}_{d}|
\label{eq}
\end{equation}
where
\begin{equation}
| A^{i}|=\left[ \frac{x-l_{i}}{q}\right]
\label{eq}
\end{equation}
\begin{equation}
| A^{i}_{d}|=\left[\frac{x-l_{i}}{dq}\right]+r_{i},d>1
\label{eq}
\end{equation}
the notation $[y]$ represents the integer part of the real number $y$. So we have
\begin{equation}
\pi(x;q,l_{i})- S_{0}(A^{i};P(z),\sqrt{x})\leq\sqrt{x}
\label{eq}
\end{equation}
namely
\begin{equation}
\pi(x;q,l_{i})= S_{0}(A^{i};P(z),\sqrt{x})+O(\sqrt{x})
\label{eq}
\end{equation}
when $l_{i}<q$,
\begin{equation}
\left[ \frac{x}{q}\right]-\left[ \frac{x-l_{i}}{q}\right]\leq 2
\label{eq}
\end{equation}
the error be very small,when $x$ be larger,it's negligible compared to $\sqrt{x}$,so we can introduce a sieve function that is very similar to $S_{0}(A^{i};P(z),\sqrt{x})$
\begin{equation}
S(A^{i};P(z),\sqrt{x})=| A|+\sum _{1<d|(P(z))}\mu (d) | A^{i}_{d}|
\label{eq}
\end{equation}
where
\begin{equation}
| A|=\left[ \frac{x}{q}\right]
\label{eq}
\end{equation}
\begin{equation}
| A^{i}_{d}|=\left[\frac{x-l_{i}}{dq}\right]+r_{i},d>1
\label{eq}
\end{equation}
so we also have
\begin{equation}
\pi(x;q,l_{i})= S(A^{i};P(z),\sqrt{x})+O(\sqrt{x})
\label{eq}
\end{equation}
According to the identity of the sieve function transformation,these any two sieve functions $S(A^{i};P(z),\sqrt{x})$, $S(A^{j};P(z),\sqrt{x})$ will be all similar sieve function $S(A^{i};P(z),\sqrt{x})\sim S(A^{j};P(z),\sqrt{x})$,namely, there must be a sieve function transformation $S(A^{i};P(z),\sqrt{x})\to S(A^{j};P(z),\sqrt{x})$,according to Theorem 2, we have
\begin{equation}
S(A^{i};P(z),\sqrt{x})-S(A^{j};P(z),\sqrt{x})\ll\sum _{p\leq z}\left||A^{i}_{p}|-|A^{j}_{p}|\right|+R_{0}
\label{eq}
\end{equation}
where
\begin{equation}
R_{0}\ll\sum _{p\leq z}\log x\ll\sqrt{x}
\label{eq}
\end{equation}
because,for $(q,l_{i})=1$, $(q,l_{j})=1$, $l_{i},l_{j}<q$,$r_{i}\ll 2$
\begin{equation}
 \sum _{p\leq \sqrt{x}}\left||A^{i}_{p}|-|A^{j}_{p}|\right|\ll\sum _{p\leq \sqrt{x}}\left|\left[\frac{x-l_{i}}{pq}\right]-\left[\frac{x-l_{j}}{pq}\right]\right|+\sum _{p\leq \sqrt{x}}2\ll\sqrt{x}
\label{eq}
\end{equation}
so we have
\begin{equation}
S(A^{i};P(z),\sqrt{x})-S(;A^{j}P(z),\sqrt{x})\ll\sqrt{x}
\label{eq}
\end{equation}
combining with (118) (122), we can get the following relations
\begin{equation}
\pi(x;q,l_{i})-\pi(x;q,l_{j})\ll\sqrt{x}
\label{eq}
\end{equation}
since the mean value of $\pi(x;q,l_{i})$ be equal to $\frac{\pi(x)}{\phi(q)}$,thus for $x>2$,we have
\begin{equation}
\pi(x;q,l_{i})=\frac{\pi(x)}{\phi(q)}+O\left(\sqrt{x}\right)
\label{eq}
\end{equation}

\begin{equation}
\pi(x;q,l_{j})=\frac{\pi(x)}{\phi(q)}+O\left(\sqrt{x}\right)
\label{eq}
\end{equation}
where the constant $O$ is a calculable constant.

To sum up,let$(q,l)=1$, $l<q$ ,then we have
\begin{equation}
\pi(x;q,l)=\frac{\pi(x)}{\phi(q)}+O\left(\sqrt{x}\right)
\label{eq}
\end{equation}

Theorem 4 be proved.
\end{proof}

\section{Conclusion}

In this paper,we introduces the basic principle of sieve function transformation,and got some of the theorems about the similar sieve function,applying the method of sieve function transformation,the prime number theorem in arithmetic sequence is improved.The prime number theorem in the arithmetic sequence is closely related to General Riemann Hypothesis, if the theorem 4 is established, then General Riemann Hypothesis may be true.The theory of sieve function transformation can also be used for many other studies on number theory.

\end{document}